\documentclass[12pt]{article}
\usepackage{amsmath,amsthm,amsfonts,amssymb,amscd}
\usepackage{graphics}
\usepackage[pdftex]{graphicx}
\def\disp{\displaystyle}

\def\Limsup{\mathop{{\rm Lim}\,{\rm sup}}}

\def\tto{\;{\lower 1pt\hbox{$\rightarrow$}}\kern -10pt
\hbox{\raise 2pt\hbox{$\rightarrow$}}\;}
\def\Hat{\widehat}
\def\hat{\widehat}

\def\Bar{\overline}
\def\ra{\rangle}
\def\la{\langle}
\def\ve{\varepsilon}
\def\B{I\!\!B}
\def\h{\hfill\Box}
\def\R{I\!\!R}
\def\IN{I\!\!N}

\def\ox{\bar{x}}

\def\ov{\bar{v}}

\def\ou{\bar{u}}

\def\gph{\mbox{\rm gph}\,}

\def\dom{\mbox{\rm dom}\,}

\def\cl*co{\mbox{\rm cl}^*\mbox{\rm co}\,}
\def\cl{\mbox{\rm cl}\,}

\def\substack#1#2{{\scriptstyle{#1}\atop\scriptstyle{#2}}}

\def\h{\hfill\triangle}

\def\O{\Omega}

\def\ph{\varphi}

\def\emp{\emptyset}
\def\st{\stackrel}
\def\oR{\Bar{\R}}
\def\lm{\lambda}
\def\gg{\gamma}

\def\hs7{\hspace*{7pt}}

\renewcommand{\theequation}{\thesection.\arabic{equation}}

\def\h{\hfill\Box}
\def\kk{\kappa}
\begin{document}
\newtheorem{Theorem}{Theorem}[section]
\newtheorem{Proposition}[Theorem]{Proposition}
\newtheorem{Remark}[Theorem]{Remark}
\newtheorem{Lemma}[Theorem]{Lemma}
\newtheorem{Corollary}[Theorem]{Corollary}
\newtheorem{Definition}[Theorem]{Definition}
\newtheorem{Example}[Theorem]{Example}
\newtheorem{Fact}[Theorem]{Fact}
\renewcommand{\theequation}{\thesection.\arabic{equation}}
\normalsize
\def\proof{
\normalfont
\medskip
{\noindent\itshape Proof.\h
space*{6pt}\ignorespaces}}
\def\endproof{$\h$\vspace*{0.1in}}
\title{\bf Coderivative Characterizations of Maximal Monotonicity for Set-Valued Mappings}
\date{}
\author{N. H. Chieu\footnote{Department of Mathematics, Vinh University, Vinh, Nghe An, Vietnam; email: nghuychieu@gmail.com. Research of this author was partially supported by the Vietnam National Foundation for Science and Technology Development (NAFOSTED) under grant 101.01-2014.56.}, G. M. Lee\footnote{Department of Applied Mathematics, Pukyong National University, Busan 608-737, Republic of Korea; email:
gmlee@pknu.ac.kr.  Research of this author was partially supported by the Basic Science Research Program through the National Research Foundation of Korea (NRF) funded by the Ministry of Education (NRF-2013R1A1A2005378)}, B. S. Mordukhovich\footnote{Department of Mathematics, Wayne State University, Detroit, MI 48202, USA; email: boris@math.wayne.edu. Research of this author was partially supported by the USA National Science Foundation under grant DMS-12092508}, T. T. A. Nghia\footnote{Department of Mathematics and Statistics, Oakland University, Rochester, MI 48309, USA; email: nttran@oakland.edu.}}
\maketitle
\begin{quote}
{\small \noindent {\bf Abstract.} This paper concerns generalized differential characterizations of maximal monotone set-valued mappings. Using advanced tools of  variational analysis, we establish coderivative criteria for maximal monotonicity of set-valued mappings, which seem to be the first infinitesimal characterizations of maximal monotonicity outside the single-valued case. We also present second-order  necessary and sufficient conditions for lower-${\cal C}^2$ functions to be convex and strongly convex. Examples are provided to illustrate the obtained results and the imposed assumptions.}

\end{quote}
\maketitle

\section{Introduction}
\setcounter{equation}{0}

The notion of maximal monotone operators appeared in the early 1960s and since that has been well recognized as a fundamental tool in the study of various aspects of partial differential equations, optimization, equilibrium theory, etc.; particularly those concerning existence and uniqueness of solutions, stability issues, convergence of numerical algorithms, and related topics; see, e.g., \cite{BC,BZ,B,BI08,FP,Min62,Phelps89,r76,rw,S98} and the references therein. We specially emphasize a crucial role of maximal monotonicity in the theory and applications of the {\em sweeping process} (``processus du rafle") introduced and investigated by Jean Jacques Moreau \cite{more71}; see more details and recent developments in the survey paper \cite{CTh}.

Since a {\em single-valued} and {\em continuous} monotone mapping is automatically maximal monotone (see, e.g., \cite[Corollary~20.25]{BC}), the maximality issue does not arise in this case. The classical criterion of monotonicity presented in \cite[Proposition~12.3]{rw} tells us that a differentiable single-valued mapping is monotone if and only if its derivative is positive-semidefinite at every point. Infinitesimal characterizations of monotonicity for possibly nondifferentiable mappings have attracted much attention in the literature. In his landmark paper \cite{Min62}, Minty established a sufficient condition for monotonicity of nondifferentiable monotone mappings by using directional derivatives. Jiang and Qi \cite[Proposition~2.3]{JiQi95} and Luc and Schaible \cite[Proposition~2.1]{LucSch96} independently proved that a locally Lipschitzian mapping defined on an open convex set of $\R^n$ is monotone if and only if its Clarke's generalized Jacobian is pointwise positive-semidefinite at every point. Replacing generalized Jacobian matrices by their approximate Jacobian matrices, Jeyakumar et al. \cite[Theorem~3.1]{JeyLucSch98} derived a sufficient condition for monotonicity of continuous single-valued mappings between finite-dimensional spaces. More recently \cite{CT}, Chieu and Trang obtained necessary and sufficient conditions for monotonicity of continuous single-valued mappings in both finite-dimensional and infinite-dimensional settings via positive-semidefiniteness of the coderivative constructions by Mordukhovich \cite{M1} (see Sections~2), extending in this way the classical characterization of monotonicity for smooth mappings.\vspace*{0.03in}

However, the major role in nonlinear and variational analysis and their applications is played not by single-valued but {\em intrinsically set-valued} maximal monotone operators that include, e.g., subdifferential mappings for lower semicontinuous proper convex functions and normal cone mappings associated with close convex sets. In particular, such set-valued mappings allow us to adequately describe variational inequalities and complementarity problems in Robinson's framework of generalized equations \cite{rob}, the aforementioned Moreau's sweeping process, etc. It is very attractive and challenging therefore to establish verifiable infinitesimal conditions (better, complete characterizations) of maximal monotonicity and related properties for set-valued mappings in finite and infinite dimensions.

To the best of our knowledge, the first result in this direction was obtained by Poliquin and Rockafellar \cite[Theorem~2.1]{PR2} who derived a {\em necessary} condition for the maximal monotonicity of set-valued mappings between finite-dimensional spaces in terms of the positive-semidefiniteness of the limiting coderivative; this condition was extended in \cite{CH,MN2} to Hilbert spaces and reversed in \cite{CT} for single-valued mappings. Note that the motivation of \cite{PR2} came from the application to {\em tilt stability} in optimization which theory has been flowering during the recent years; see, e.g., \cite{BS,DL,DMN,EW,MN2,MN1} and the references therein.

Another impact to the study and of monotonicity properties for set-valued mappings has been recently done by Mordukhovich and Nghia \cite{MN} who established complete coderivative characterizations of {\em strong local} maximal monotonicity in finite and infinite dimensions with applications to {\em full stability} (in the Lipschitzian and H\"olderian frameworks) of parametric variational systems. The approach of \cite{MN} made used, along with advanced tools of variational analysis and generalized differentiation, {\em hypomonotonicity} properties of set-valued mappings, which will be exploited in what follows.\vspace*{0.03in}

The main goal of this paper is to establish complete {\em coderivative characterizations} of {\em maximal monotonicity} of set-valued mappings in Hilbert spaces via pointwise {\em positive-semidefiniteness} conditions and appropriate properties of global and (semi)local hypomonotonicity. The results obtained seem to be the first infinitesimal characterizations of maximal monotonicity outside the single-valued setting even in the case of finite dimensions. As consequences of these characterizations, we derive second-order necessary and sufficient conditions for lower-${\cal C}^2$ functions to be {\em convex} and {\em strongly convex}. These conditions are expressed in terms of {\em second-order subdifferenntials/generalized Hessians} of extended-real-valued functions and extend the classical result of real analysis saying that a ${\cal C}^2$ function is convex if and only if its Hessian is positive-semidefinite at any point.\vspace*{0.03in}

The paper is organized as follows. Section~2 recalls some basic notions and facts from variational analysis that are employed in the sequel. Section~3 is the main part of our analysis, which contains several coderivative characterizations of maximal monotonicity for set-valued mappings in Hilbert spaces. Section~4 is devoted to the study of convexity and strong convexity for lower-${\cal C}^2$ functions via second-order subdifferentials. Finally, we present some concluding remarks and formulate open questions in Section~5.

Our notation is standard in variational analysis and generalized differentiation; cf.\ \cite{M1,rw}. Throughout the paper we assume that $X$ is a {\em Hilbert} space being identified with its dual space. As usual, the symbol $\la\cdot,\cdot,\ra$ signifies the canonical pairing in $X$ with the norm $\|x\|=\sqrt{\la x,x\ra}$. We denote by $\B$ the closed unit ball in $X$ and by $\B_r(\ox):=\ox+r\B$ the ball with radius $r>0$ and center $\ox\in X$. The notation $\st{w}\rightarrow$ indicates for the weak convergence in $X$. Given a set-valued mapping $F\colon X\tto X$ and a point $\ou\in X$, the symbol
\begin{equation}\label{pk}
\begin{array}{ll}
\disp\Limsup_{u\to\ou}F(u):=\Big\{v\in X\Big|&\exists\;\mbox{ sequences }\;u_k\to\ou,\;v_k\st{w}{\to}v\;\mbox{ such that}\\
&v_k\in F(u_k)\;\mbox{ for all }\;k\in\IN:=\{1,2,\ldots\}\Big\}
\end{array}
\end{equation}
stands for the {\em sequential Painlev\'e-Kuratowski outer/upper limit} of $F(u)$ as $u\to\ou$.

\section{Preliminaries}
\setcounter{equation}{0}

Here we mainly follow the book \cite{M1} referring the reader also to \cite{BZ,rw} for related and additional material. Given a proper (i.e., not identically equal to infinity) extended-real-valued function $f\colon X\to\oR:=\R\cup\{\infty\}$ and its domain point $\ou\in\dom f:=\{u\in X|\;f(u)<\infty\}$, the {\em regular/Fr\'echet subdifferential} (known also as the presubdifferential and the viscosity subdifferential) of $f$ at $\ou$ is
\begin{eqnarray}\label{2.1}
\Hat\partial f(\ou):=\Big\{v\in X\Big|\;\liminf_{u\to\ou}\frac{f(u)-f(\ou)-\la v,u-\ou\ra}{\|u-\ou\|}\ge 0\Big\}
\end{eqnarray}
with $\hat\partial f(\ou):=\emp$ if $\ou\notin\dom f$. The {\em limiting/Mordukhovich subdifferential} (known also as the basic subdifferential) of $f$ at $\ou\in\dom f$ is defined via \eqref{pk} by
\begin{equation}\label{2.2}
\partial f(\ou):=\Limsup_{u\st{f}\to\ou}\hat\partial f(u),
\end{equation}
where the notation $u\st{f}{\to}\ou$ means that $u\to\ou$ with $f(u)\to f(\ou)$.  It is well known that both regular and limiting subdifferential reduce to the classical subdifferential of convex analysis when the function $f$ is convex. On the other hand, the limiting subdifferential of nonconvex functions and related normal cone/coderivative constructions for sets and mappings enjoy {\em full calculus}, which is not the case for (\ref{2.1}) and its set/mapping counterparts.

Given further a set $\O\subset X$ with its indicator function $\delta(u;\O)$ equal to $0$ for $u\in\O$ and to $\infty$ otherwise, the regular and limiting {\em normal cones} to $\O$ at $\ou\in\O$ are defined, respectively, via the corresponding subdifferentials \eqref{2.1} and \eqref{2.2} by
\begin{eqnarray}\label{2.3}
\Hat N(\ou;\O):=\Hat\partial\delta(\ou;\O)\;\mbox{ and }\;N(\ou;\O):=\partial\delta(\ou;\O).
\end{eqnarray}

Now we consider a set-valued mapping $F:X\tto X$ and associate with it the {\em domain} $\dom F$ and the {\em graph} $\gph F$  by
$$
\dom F:=\big\{u\in X\big|\;F(u)\ne\emp\big\}\;\mbox{ and }\;\gph F:=\big\{(u,v)\in X\times X\big|\;v\in F(u)\big\}.
$$
The mapping $F$ is said to be {\em proper} when $\dom F\ne\emp$, which is always assumed. Define by
\begin{eqnarray}\label{2.4}
\Hat D^*F(\ou,\ov)(w):=\big\{z\in X|\;(z,-w)\in\Hat N((\ou,\ov);\gph F)\big\}\;\mbox{ for all }\;w\in X
\end{eqnarray}
the {\em regular coderivative} of $F$ at $(\ou,\ov)\in\gph F$ and by
\begin{eqnarray}\label{2.5}
D_M^*F(\ou,\ov)(w):=\Limsup_{\substack{(u,v)\to(\ou,\ov)}{y\to w}}\Hat D^* F(u,v)(y)\;\mbox{ for all }\;w\in X
\end{eqnarray}
the {\em mixed limiting coderivative} of $F$ at $(\ou,\ov)$, where the convergence $y\to w$ is strong in $X$ while the outer limit in \eqref{2.5} is taken by \eqref{pk} in the weak topology of $X$; see \cite{M1} for more discussions. We omit the subscript $``M"$ in \eqref{2.5} when $X$ is finite-dimensional and also drop indicating $\ov=F(\ou)$ if $F$ is single-valued. When $F$ is single-valued and continuously differentiable around $\ou$ (or strictly differentiable at this point), we get
\begin{eqnarray*}
\Hat D^*F(\ou)(w)=D^*_M F(\ou)(w)=\big\{\nabla F(\ou)^*w\big\}\;\mbox{ for all }\;w\in X
\end{eqnarray*}
via the adjoint derivative operator $\nabla F(\ou)^*$; see, e.g., \cite[Theorem~1.38]{M1}.

Next we recall two second-order subdifferential/generalized Hessian constructions introduced by the scheme suggested in \cite{mor92} as a coderivative of a first-order subdifferential mapping; see \cite{M1,MN2,MR} for more details and discussions.

\begin{Definition} [\bf second-order subdifferentials]\label{2nd} Let $f:X\to\oR$ with $\ou\in\dom f$, and let $\ov\in\partial f(\ou)$. Then we say that:

{\bf(i)} The {\sc combined second-order subdifferential} of $f$ at $\ou$ relative to $\ov$ is the set-valued mapping $\breve\partial^2 f(\ox,\ov): X\tto X$ with the values
\begin{equation}\label{2.6}
\breve\partial^2f(\ou,\ov)(w):=\big(\Hat D^*\partial f\big)(\ou,\ov)(w)\;\mbox{ for all }\;w\in X.
\end{equation}

{\bf (ii)} The {\sc mixed second-order subdifferential} of $f$ at $\ou$ relative to $\ov$ is the set-valued mapping $\partial^2_M f(\ou,\ov): X\tto X$ with the values
\begin{eqnarray}\label{2.7}
\partial_M^2f(\ou,\ov)(w):=\big(D_M^*\partial f\big)(\ou,\ov)(w)\;\mbox{ for all }\;w\in X.
\end{eqnarray}
\end{Definition}

It is worth mentioning that if $f$ is ${\cal C}^2$ around $\ou$ with $\ov=\nabla f(\ou)$, then both $\breve\partial^2f(\ou,\ov)(w)$ and $\partial_M^2f(\ou,\ov)(w)$ reduce to the classical (symmetric) single-valued Hessian operator:
$$
\breve\partial^2f(\ou,\ov)(w)=\partial_M^2 f(\ou,\ov)(w)=\big\{\nabla^2f(\ou)^*w\big\}=\big\{\nabla^2f(\ou)w\big\}\;\mbox{ for all }\; w\in X.
$$

One of the main fact of generalized differentiation largely employed in our paper is the following {\em mean value inequality} for Lipschitz continuous functions; \cite[Corollary~3.50(ii)]{M1}. For the reader's convenience we formulate it here.\\[1ex]
{\bf Mean-value inequality.} {\em Let $f:X\to\oR$ be a lower semicontinuous $($l.s.c.$)$ function with $a\in\dom f$. Then for any $b\in X$ and $\ve>0$  we have the estimate
\begin{equation}\label{mvt}
|f(b)-f(a)|\le\|b-a\|\sup\big\{\|v\|\;\big|\;v\in\hat\partial f(c),\;c\in[a,b]+\ve\B\big\},
\end{equation}
where $[a,b]:=\{\lm a+(1-\lm)b|\;\lm\in[0,1]\}$.}

\section{Characterizations of Maximal Monotonicity}
\setcounter{equation}{0}

The following notions of (global) monotonicity for set-valued mappings are main objects of our study and applications in this paper.

\begin{Definition}[\bf monotone set-valued operators]\label{Mono} Given $T:X\tto X$, we say that:

{\bf (i)} $T$ is {\sc monotone} on $X$ if
\begin{eqnarray}\label{3.1}
\la v_1-v_2,u_1-u_2\ra\ge 0\;\mbox{ for all }\;(u_1,v_1),\;(u_2,v_2)\in \gph T.
\end{eqnarray}
$T$ is said to be {\sc maximal monotone} on $X$ if in addition we have $\gph T=\gph S$ whenever $S$ is monotone with $\gph T\subset\gph S$.

{\bf (ii)} $T$ is {\sc  hypomonotone} on $X$ if there exists a number $r>0$ such that $T+r I$, where $I\colon X\to X$ is the identity mapping. This means that
\begin{eqnarray}\label{3.2}
\la v_1-v_2,u_1-u_2\ra\ge-r \|u_1-u_2\|^2\;\mbox{ for all }\;(u_1,v_1),(u_2,v_2)\in\gph T.
\end{eqnarray}
\end{Definition}

First we present a characterization of maximal monotonicity via the positive-semidefiniteness condition for the regular coderivative and global hypomonotonicity.

\begin{Theorem}[\bf regular coderivative and global hypomonotonicity characterization of maximal monotonicity]\label{thm1} Let $T:X\tto X$ be a set-valued mapping with closed graph in the norm topology of $X\times X$. The following assertions are equivalent:

{\bf (i)} $T$ is maximal monotone on $X$.

{\bf (ii)} $T$ is hypomonotone on $X$ and for any $(u,v)\in\gph T$ we have
\begin{eqnarray}\label{Cod}
\la z,w\ra \ge 0\;\mbox{ whenever }\;z\in\Hat D^*T(u,v)(w).
\end{eqnarray}
\end{Theorem}
\noindent{\bf Proof.} Since the monotonicity obviously yields hypomonotonicity, implication (i)$\Longrightarrow$(ii) follows from  \cite[Lemma~5.2]{CH} and also from \cite[Lemma~6.2]{MN2}).

To verify the converse implication (ii)$\Longrightarrow$(i), suppose that $T$ is hypomonotone and that condition \eqref{Cod} is satisfied. Then there is some number $r>0$ such that $T+rI$ is monotone. Take any $s>r$ and define $F:X\tto X$ by $\gph F:=\gph (T+sI)^{-1}$. For any $(v_i,u_i)\in\gph F$, $i=1,2$ we have $(u_i,v_i-su_i)\in\gph T$ and thus deduce from \eqref{3.2} that
\[
\la v_1-s u_1-v_2+su_2,u_1-u_2\ge -r\|u_1-u_2\|^2.
\]
The latter implies in turn that the inequalities
\begin{eqnarray*}
\|v_1-v_2\|\cdot\|u_1-u_2\|\ge\la v_1-v_2,u_1-u_2\ra\ge(s-r)\|u_1-u_2\|^2,
\end{eqnarray*}
which allow us to arrive at the estimate
\begin{eqnarray}\label{Lip}
\|u_1-u_2\|\le\frac{1}{s-r}\|v_1-v_2\|
\end{eqnarray}
verifying that $F$ is single-valued and Lipschitz continuous with modulus $(s-r)^{-1}$ on its domain. To proceed further, fix any $z\in X$ and define $f_z:X\to\oR$ by
\begin{eqnarray}\label{f}
f_z(v):=\left\{\begin{array}{ll}
\la z, F(v)\ra&\mbox{if }\;v\in\dom F,\\
\infty&\mbox{otherwise}.
\end{array}\right.
\end{eqnarray}
Since $\gph T$ is closed, it is easy to check that $\gph F$ is also closed in $X\times X$. Next we show that $f_z$ is lower semicontinuous on $X$. Arguing by contradiction, suppose that there exist $\ve>0$ and a sequence $v_k$ converging to some $v\in X$ such that $f_z(v_k)<f_z(v)-\ve$. If $f_z(v)=\infty$, then $v\notin\dom F$ and $v_k\in\dom F$. It follows from \eqref{Lip} that $\|F(v_k)-F(v_j)\|\le(s-r)^{-1}\|v_k-v_j\|$, and so $\{F(v_k)\}$ is a Cauchy sequence converging to some $u\in X$. Hence the sequence $(v_k,F(v_k))\in\gph F$ converges to $(v,u)\in\gph F$ due to the closedness of $\gph F$. This gives us $F(v)=u$ and contradicts $v\notin\dom F$. In the remaining case of $f_z(v)<\infty$ we get from \eqref{Lip} and \eqref{f} the estimates
\[
|f_z(v_k)-f_z(v)|\le\|z\|\cdot\|F(v_k)-F(v)\|\le\|z\|\cdot\frac{1}{s-r}\|v_k-v\|\to 0\;\mbox{ as }\;k\to\infty,
\]
which is also a contradiction due to the assumption $f_z(v_k)<f_z(v)-\ve$. This justifies the lower semicontinuity of $f_z$ on the space $X$ for any fixed $z\in X$.

Now we claim that $T$ is {\em monotone}. To proceed, pick two pair $(u_i,v_i)\in\gph T$ and get
$$
(y_i,u_i)\in\gph F\;\mbox{ with }\;y_i:=v_i+su_i,\;i=1,2.
$$
Applying the mean value inequality (\ref{mvt}) to the l.s.c.\ function $f_z$ tells us that
\begin{equation}\label{Lip2}
|\la z,u_1-u_2\ra|=|f_z(y_1)-f_z(y_2)|\le\|y_1-y_2\|\sup\big\{\|w\|\;\big|\;w\in\hat\partial f_z(y),\;y\in [y_1,y_2]+\ve\B\big\}
\end{equation}
with $[y_1,y_2]:=\{\lm y_1+(10-\lm)y_2|\;\lm\in [0,1]\}$ and fixed $\ve>0$. Since $\hat\partial f_z(y)=\emp$ if $y\notin\dom f_z$, it suffices to consider the case of $y\in\dom f_z\cap([y_1,y_2]+\ve\B)=\dom F\cap([y_1,y_2]+\ve\B)$ in \eqref{Lip2}. Take any $y$ from the latter set and observe that
\begin{equation}\label{incl}
w \in\hat D^*F(y)(z)\;\mbox{ whenever }\;w\in\hat\partial f_z(y).
\end{equation}
Indeed, it follows from the definition of the regular subgradient $w\in\hat\partial f_z(y)$ that
$$
\liminf\limits_{v\to y}\frac{f_z(v)-f_z(y)-\la w,v-y\ra}{\|v-y\|}\ge 0,
$$
which can be equivalently written by the construction of $f_z$ in (\ref{f}) as
$$
\disp\liminf_{v\st{{\rm{dom}}\,F}{\to}y}\frac{\la z,F(v)\ra-\la z,F(y)\ra-\la w,v-y\ra}{\|v-y\|}\ge 0.
$$
The latter readily implies that
$$
\disp\liminf_{(v,u)\st{{\rm{gph}}\,F}{\to}(y,F(y))}\frac{\la z,u- F(y)\ra-\la w,v-y\ra}{\|v-y\|+\|u-F(y)\|}\ge 0.
$$
Hence we get from \eqref{2.1}, \eqref{2.3}, and \eqref{2.4} that
$$
(w,-z)\in\Hat N\big((y,F(y));\gph F\big)\Longleftrightarrow w\in\hat D^*F(y)(z)=\hat D^*(T+sI)^{-1}(y)(z)
$$
and therefore $-z\in\hat D^*(T+sI)(F(y),y)(-w)$. It easily follows from the coderivative sum rule in \cite[Theorem~1.62]{M1} that
\begin{equation}\label{co}
-z+sw\in\hat D^*T\big(F(y),y-sF(y)\big)(-w).
\end{equation}
Combining this with \eqref{Cod} tells us that $\la-z+sw,-w\ra\ge 0$, which yields
\begin{eqnarray}\label{ZW}
\|z\|\cdot\|w\|\ge\la z,w\ra\ge s\|w\|^2
\end{eqnarray}
and implies furthermore together with the estimate \eqref{Lip2} that
\begin{eqnarray*}
|\la z, u_1-u_2\ra|\le s^{-1}\|z\|\cdot\|y_1-y_2\|.
\end{eqnarray*}
Since this inequality holds for all $z\in X$, we get
\begin{eqnarray*}
\|u_1-u_2\|\le s^{-1}\|y_1-y_2\|=s^{-1}\|v_1+su_1-v_2-su_2\|
\end{eqnarray*}
and then deduce by the elementary transformation that
\begin{eqnarray*}
s^2\|u_1-u_2\|\le\|(v_1-v_2)+s(u_1-u_2)\|^2=\|v_1-v_2\|^2+2\la v_1-v_2,u_1-u_2\ra+s^2\|u_1-u_2\|^2.
\end{eqnarray*}
Therefore we arrive at the inequality
\[
0\le\frac{1}{2s}\|v_1-v_2\|^2+\la v_1-v_2,u_1-u_2\ra\;\mbox{ for any }\;s>r.
\]
Letting there $s\to\infty$ shows that
\begin{equation*}\label{mon}
0\le\la v_1-v_2,u_1-u_2\ra\;\mbox{ for all }\;(u_1,v_1),(u_2,v_2)\in\gph T
\end{equation*}
and thus verifies the monotonicity of $T$.

It remains to prove that $T$ is {\em maximal monotone}. Since $T$ is proper, there exists a pair $(u_0,v_0)\in\gph T$ such that
$$
u_0=(T+sI)^{-1}(y_0)\;\mbox{ with }\;y_0:=v_0+s u_0.
$$
Applying again the mean value inequality (\ref{mvt}) to the function $f_z$ from (\ref{f}), we have
\begin{eqnarray}\label{Lip4}
|f_z(y)-f_z(y_0)|\le\|y-y_0\|\sup\big\{\|w\|\;\big|\;w\in\hat\partial f_z(x),\;x\in[y,y_0]+\ve\B\big\}
\end{eqnarray}
for any $y\in X$. It follows similarly to \eqref{ZW} that $\|w\|\le s^{-1}\|z\|$ for all $w\in\hat\partial f_z(x)$ with $x\in\dom F\cap([y,y_0]+\ve\B)$. This together with \eqref{Lip4} gives us the estimates
\[
|f_z(y)-f_z(y_0)|\le s^{-1}\|z\|\cdot\|y-y_0\|.
\]
Hence $f_z(y)<\infty$ and so $F(y)\ne\emp$ for all $y\in X$, which means that $\dom(T+sI)^{-1}=X$. Employing now the classical Minty theorem (see, e.g., \cite[Theorem~4.4.7 and Remark~4.4.8]{BI08}) and taking into account the monotonicity of $T$ justified above, we conclude that $T$ is maximal monotone and thus complete the proof of the theorem.\endproof

Our next goal is to obtain another version of the coderivative characterization in Theorem~\ref{thm1} with replacing the global hypomonotonicity of $T$ in assertion (ii) therein by a certain local hypomonotonicity. Besides being interesting for its own sake, it is needed for the subsequent applications in Section~4 to characterize convexity and strong convexity of lower-${\cal C}^2$ functions. In fact, for these purposes we need to modify the conventional notion of local monotonicity and hypomonotonicity, which are dealing with neighborhoods in the {\em product} space $X\times X$; see, e.g., \cite{MN,P} with more references and discussions. The notions we use in what follows concern neighborhoods only in the {\em domain} space $X$. Such a local monotonicity has been considered in \cite[Example~12.28]{rw}. We will name the domain versions as {\em semilocal} monotonicity and hypomonotonicity, which reflects their nature and distinguishes them from their fully localized product counterparts.

\begin{Definition}[\bf semilocal monotonicity and hypomonotonicity]\label{LH} We say that the mapping $T\colon X\tto X$ is {\sc semilocally hypomonotone} $($resp.\ {\sc semilocally monotone}$)$ at $\bar u\in\dom T$ if there exist a neighborhood $U$ of $\ou$ and a number $r>0$ $($resp.\ $r=0)$ with
\begin{equation}\label{lh}
\la v_1-v_2,u_1-u_2\ra\ge-r\|u_1-u_2\|^2\;\mbox{ for all }\;(u_1,v_1),(u_2,v_2)\in\gph T\cap(U\times X).
\end{equation}
Given a set $\O\subset X$, we say that $T$ is semilocally hypomonotone $($resp.\ monotone$)$ {\sc on} $\O$ if it is semilocally hypomonotone $($resp.\ monotone$)$ at every point $\bar u\in\O\cap\dom T$.
\end{Definition}

Establishing the desired semilocal version of Theorem~\ref{thm1} requires an additional convexity assumption on the domain of $T$, which is shown below to be essential by providing a counterexample. To proceed in this direction, we first present the following lemma proved in \cite[Theorem~5]{Min62} by similar arguments for single-valued operators.

\begin{Lemma}[\bf semilocal monotonicity of set-valued mappings with convex domains]\label{lmMainthm} Let $T\colon X\tto X$ be a semilocally monotone mapping on $X$, and let its domain $\dom T$ be convex. Then $T$ is $($globally$)$ monotone on $X$.
\end{Lemma}
\noindent{\bf Proof.} Pick any $(u_1,v_1),(u_2,v_2)\in\gph T$ and get $[u_1,u_2]\subset\dom T$ by the convexity of $\dom T$. Since $T$ is semilocally monotone, for each $x\in[u_1,u_2]$ there is $\gamma_x>0$ such that
\begin{equation}\label{cover}
\la y_1-y_2,x_1-x_2\ra\ge 0\;\mbox{ whenever }\;(x_1,y_1),(x_2,y_2)\in\gph T\cap\big({\rm int}\,\B_{\gg_x}(x)\times X\big).
\end{equation}
The compactness of $[u_1,u_2]$ allows us to select $x_i\in [u_1,u_2],$ $i=1,\ldots,n$ satisfying
$$
[u_1,u_2]\subset\bigcup\limits_{i=1}^n\Big({\rm int}\,\B_{\gamma_{x_i}}(x_i)\Big).
$$
Thus we can find $0=t_0<t_1<\ldots<t_m=1$ such that for each $j\in\{0,\ldots,m-1\}$ it holds
$$
[\hat u_j,\hat u_{j+1}]\subset{\rm int}\,\B_{\gamma_{x_{i_j}}}(x_{i_j})\;\mbox{ with some }\;i_j\in\big\{1,\ldots,n\big\},
$$
where $\hat u_j:=u_1+t_j(u_2-u_1)$. Since $\hat u_j\in[u_1,u_2]\subset\dom T$ for each $j\in\{0,\ldots,m\}$, there exist $\hat v_j\in T(\hat u_j)$  with $\hat v_0=v_1$ and $\hat v_m=v_2$. It follows from \eqref{cover} that
$$
(t_{j+1}-t_j)\la\hat v_{j+1}-\hat v_j,u_2-u_1\ra=\langle\hat v_{j+1}-\hat v_j,\hat u_{j+1}-\hat u_j\rangle\ge 0,
$$
which implies that $\la \hat v_{j+1}-\hat v_j,u_2-u_1\ra\ge 0$ whenever $j\in\{0,\ldots,m-1\}$. Hence we get
$$
\la v_2-v_1,u_2-u_1\ra=\sum\limits_{j=0}^{m-1}\la\hat v_{j+1}-\hat v_j,u_2-u_1\ra\ge 0
$$
and thus verify the global monotonicity of the operator $T$.\endproof

Now we are ready to obtain a semilocal counterpart of the coderivative characterization in Theorem~\ref{thm1} under the convexity assumption on $\dom T$. Example~\ref{ex} below demonstrates that the latter assumption cannot be dropped. Since the proof of the following theorem is similar in some places to that of Theorem~\ref{thm1}, we omit the corresponding details.

\begin{Theorem}[\bf regular coderivative and semilocal hypomonotonicity characterization of maximal monotonicity]\label{thm2} Let $T:X\tto X$ be a set-valued mapping with closed graph and convex domain. The following assertions are equivalent:

{\bf (i)} $T$ is maximal monotone on $X$.

{\bf (ii)} $T$ is semilocally  hypomonotone on $X$ and satisfies the regular coderivative condition \eqref{Cod} for any $(u,v)\in\gph T$.
\end{Theorem}
\noindent{\bf Proof.} Implication (i)$\Longrightarrow$(ii) follows from Theorem~\ref{thm1}. To verify the converse implication, suppose that condition \eqref{Cod} holds and that $T$ is semilocally hypomonotone. The latter allows us to find, for each $\ou\in\dom T$, positive numbers $\delta$ and $r$ such that
\begin{equation}\label{hy}
\la v_1-v_2,u_1-u_2\ra\ge-r\|u_1-u_2\|^2\;\mbox{ whenever }\;(u_1,v_1),(u_2,v_2)\in\gph T\cap(\B_\delta(\ou)\times X).
\end{equation}
Take any $s>r$ and define the mapping $F:X\tto X$ by $\gph F:=\gph(T+sI)^{-1}\cap(X\times\B_\delta(\ou))$. Picking arbitrarily pairs $(v_i,u_i)\in\gph F$, $i=1,2$ we have $(u_i,v_i-su_i)\in\gph T\cap(\B_\delta(\ou)\times X)$. It follows from \eqref{hy} that
\[
\la v_1-su_1-v_2+su_2,u_1-u_2\ra\ge-r\|u_1-u_2\|^2.
\]
Similarly to \eqref{Lip2} we deduce from the latter that
\begin{equation}\label{Lip5}
\|u_1-u_2\|\le\frac{1}{s-r}\|v_1-v_2\|\;\mbox{ for all }\;(v_1,u_1),(v_2,u_2)\in\gph F.
\end{equation}
This implies that $F$ is single-valued and Lipschitz continuous on $\dom F$. For any fixed vector $z\in X$ we also define the function $f_z:X\to\oR$ as in \eqref{f} and prove similarly to Theorem~\ref{thm1} that $f_z$ is lower semicontinuous on $X$.

Now pick arbitrary pairs $(u_1,v_1),(u_2,v_2)\in\gph T\cap({{\rm int}\,\B_\delta(\ou)}\times X)$ and fix $\ov\in T(\ou)$.  Then $F(y_i)=u_i\in\B_\delta(\ou)$ with $y_i:=v_i+s u_i$.  Applying the mean value inequality (\ref{mvt}) for any $\ve\in(0,\sqrt{s})$ tells us that
\begin{equation}\label{Lip6}
|\la z,u_1-u_2\ra=|f_z(y_1)-f_z(y_2)|\le\|y_1-y_2\|\sup\big\{\|w\|\big|\;w\in\hat\partial f_z(y),y\in [y_1,y_2]+\ve\B\big\}.
\end{equation}
Similar to \eqref{incl} we get $\hat\partial f_z(y)\subset\hat D^*F(y)(z)$ for all $y\in\dom F\cap([y_1,y_2]+\ve\B)$ and then for any $y\in\dom F\cap([y_1,y_2]+\ve\B)$ find some $y_0\in\ve\B$ and $t\in[0,1]$ satisfying $y=ty_1+(1-t)y_2+y_0$. Since $F(\ov+s\ou)=\ou$, it follows from \eqref{Lip5} that
\begin{eqnarray}
\begin{array}{ll}\label{Lip7}
&\|F(y)-\ou\|\disp=\|F(ty_1+(1-t)y_2+y_0)-F(\ov+s\ou)\|\\
&\disp\le \frac{1}{s-r}\|ty_1+(1-t)y_2+y_0-\ov-s\ou\|\\
&=\disp\frac{1}{s-r}\|t(v_1+su_1)+(1-t)(v_2+su_2)+y_0-\ov-s\ou\|\\
&=\disp\frac{1}{s-r}\|t(v_1-\ov)+st(u_1-\ou))+(1-t)(v_2-\ov)+s(1-t)(u_2-\ou)+y_0\|\\
&\disp\le\frac{1}{s-r}\Big[t\|v_1-\ov\|+(1-t)\|v_2-\ov\|+st\|u_1-\ou\|+s(1-t)\|u_2-\ou\|+\|y_0\|\Big]\\
&\disp\le\frac{1}{s-r}\Big[\max\big\{\|v_1-\ov\|,\|v_2-\ov\|\big\}+\ve\Big]+\frac{s}{s-r}\max\big\{\|u_1-\ou\|,\|u_2-\ou\|\big\}\\
&\disp\le\frac{1}{s-r}\Big[\max\big\{\|v_1-\ov\|,\|v_2-\ov\|\big\}+\sqrt{s}\Big]+\frac{s}{s-r}\max\big\{\|u_1-\ou\|,\|u_2-\ou\|\big\}.
\end{array}
\end{eqnarray}
Since  the choice of $(u_1,v_1),(u_2,v_2),(\ou,\ov)\in\gph T\cap({{\rm int}\,\B_\delta(\ou)}\times X)$ was independent of the parameter $s>r$ and by $\max\{\|u_1-\ou\|,\|u_2-\ou\|\}<\delta$, we can find $M$ so large that
\[
\frac{1}{s-r}\max\big\{\|v_1-\ov\|,\|v_2-\ov\|+\sqrt{s}\big\}+\frac{s}{s-r}\max\big\{\|u_1-\ox\|,\|u_2-\ox\|\big\}<\delta\;\mbox{ if }\;s>M.
\]
This together with \eqref{Lip7} ensures that $F(y)\in{\rm int}\,\B_\delta(\bar x)$ and thus
\[
\Hat N\big((y,F(y));\gph F\big)=\Hat N\big((y,F(y));\gph (T+sI)^{-1}\cap(X\times \B_\delta(\ox)\big)=\Hat N\big((y,F(y));\gph (T+sI)^{-1}\big),
\]
which clearly implies in turn the equality
\[
\Hat D^*F(y)(z)=\Hat D^*(T+sI)^{-1}\big(y,F(y)\big)(z).
\]
Similarly to \eqref{co}, for any $w\in\hat\partial f_z(y)\subset\Hat D^*F(y)(z)$ we get from the latter that $-z+sw\in\Hat D^* T(F(y),y-sF(y))(-w)$.
It follows from \eqref{Cod} that $\la-z+sw,-w\ra\ge 0$, which yields
$$\|z\|\cdot\|w\|\ge\la z,w\ra\ge s\|w\|^2,\;\mbox{ i.e., }\;\|z\|\ge s\|w\|.
$$
This together with \eqref{Lip6} tells us that
\[
\la z,u_1-u_2\ra \le s^{-1}\|y_1-y_2\|\cdot\|z\|.
\]
Since the latter holds for any $z\in X$, we have
\[
\|u_1-u_2\|^2\le s^{-2}\|y_1-y_2\|=s^{-2}\|v_1+su_1-v_2-su_2\|^2=s^{-2}\|(v_1-v_2)+s(u_1-u_2)\|^2
\]
and hence arrive at the estimate
\[
0\le\frac{1}{s}\|v_1-v_2\|^2+2\la v_1-v_2,u_1-u_2\ra\;\mbox{ when }\;s>M.
\]
Letting there $s\to\infty$ shows that
\[
0\le\la v_1-v_2,u_1-u_2\ra\;\mbox{ for all }\;(u_1,v_1),(u_2,v_2)\in\gph T\cap\big({\rm int}\,\B_\delta(\ou)\times X\big),
\]
which verifies the semilocal monotonicity of $T$ at any $\ou\in\dom T$. Since $\dom T$ is convex, Lemma~\ref{lmMainthm} tells us that $T$ is globally monotone. Now we are in a position to apply Theorem~\ref{thm1} and conclude therefore that $T$ is maximal monotone on $X$.\endproof

It is well known in monotone operator theory that the maximal monotonicity of $T$ always yields the convexity of the {\em closure} of the domain ${\rm cl}(\dom T)$; see, e.g., \cite[Corollary~21.12]{BC}. This naturally gives a raise to the question: whether Theorem~\ref{thm2} is true when the condition on the convexity of $\dom T$ is replaced by the convexity of ${\rm cl}(\dom T)$? The following simple example shows that it is not true and consequently that the convexity assumption on $\dom T$ in Theorem~\ref{thm2} cannot be dropped.

\begin{Example}[\bf semilocal monotonicity does not yield the convexity of the domain]\label{ex}{\rm Define the mapping $T:\R\tto\R$ by
$$
T(x):=\begin{cases}\Big\{-\disp\frac{1}{x}\Big\}&\mbox{if }\;x\in\R\backslash\{0\},\\
\emp&\mbox{if }\;x=0.\end{cases}
$$
Observe that $\gph T$ is closed, $T$ is semilocally monotone on $\R$, $\dom T=\R\backslash\{0\}$ is nonconvex while ${\rm cl}(\dom T)=\R$ is convex. Moreover, it is obvious that assertion (ii) of Theorem~\ref{thm2} is valid, but $T$ is not globally monotone on $\R$.}
\end{Example}

The next theorem provides other coderivative characterizations of maximal monotonicity, where the regular coderivative condition \eqref{Cod} is replaced by the positive-semidefiniteness conditions imposed on the {\em mixed limiting coderivative} \eqref{2.5}. These characterizations are clearly equivalent to those presented in Theorems~\ref{thm1} and \ref{thm2}, but in this paper is more convenient for us to derive them by passing to the limit in (\ref{Cod}). Note that the limiting coderivative characterizations have a strong advantage in comparison with (\ref{Cod}) due to well-developed calculus rules for (\ref{2.5}); see Remark~\ref{rem} and Section~5 for more discussions.

\begin{Theorem}[\bf limiting coderivative characterizations of maximal monotonicity]\label{Coro1} Let $T:X\tto X$ be a set-valued mapping with closed graph. The following are equivalent:

{\bf (i)} $T$ is maximal monotone on $X$.

{\bf (ii)} $T$ is hypomonotone on $X$ and for any $(u,v)\in\gph T$ we have
\begin{eqnarray}\label{Lim-Cod}
\la z,w\ra\ge 0\;\mbox{ whenever }\;z\in D_M^*T(u,v)(w),\;w\in X.
\end{eqnarray}
If in addition the operator domain $\dom T$ is convex, then the $($global$)$ hypomonotonicity in assertion {\bf(ii)} can be equivalently replaced by the semilocal one.
\end{Theorem}
\noindent{\bf Proof.} Implication (ii)$\Longrightarrow$(i) is straightforward from Theorem~\ref{thm1} due to
\[
\Hat D^*T(u,v)(w)\subset D_M^*T(u,v)(w)\;\mbox{ for all }\;(u,v)\in\gph T,\;w\in X.
\]
Thus \eqref{Cod} follows from \eqref{Lim-Cod}, and $T$ is maximal monotone by Theorem~\ref{thm1}.

To justify the reverse implication (i)$\Longrightarrow$(ii), suppose that (i) holds, and so \eqref{Cod} is valid due to Theorem~\ref{thm1}. Picking any $(u,v)\in\gph T$ and $z\in D^*_M T(u,v)(w)$ and using definition (\ref{2.5}) of the mixed limiting coderivative, we find sequences $(u_k,v_k)\st{{\rm gph}\,T}\to(u,v)$ with $z_k\st{w}\to z$ and $w_k\to w$ satisfying $z_k\in\Hat D^*T(u_k,v_k)(w_k)$ for all $k\in\IN$. It follows from \eqref{Cod} that $\la z_k,w_k\ra\ge 0$. Letting $k\to\infty$ and taking into account that sequence $\{w_k\}$ converges {\em strongly} in $X$ give us that $\la z,w\ra\ge 0$, which verifies \eqref{Lim-Cod}.

If $\dom T$ is convex, we proceed in the same way with replacing hypomonotonicity by semilocal hypomonotonicity and using Theorem~\ref{thm2} instead of Theorem~\ref{thm1}.
\endproof

\begin{Remark}[\bf advantages of limiting coderivative characterizations]\label{rem} {\rm Although the coderivative conditions (\ref{Cod}) and (\ref{Lim-Cod}) married with the corresponding hypomonotonicity give us the equivalent characterizations of maximal monotonicity, the limiting coderivative condition (\ref{Lim-Cod}) has clear advantages in comparison with the regular coderivative one (\ref{Cod}). This is due to the well-developed {\em full calculus} for the limiting coderivative (in contrast to its regular counterpart) presented in the first volume of \cite{M1}.
The comprehensive calculus rules developed for (\ref{2.5}) allow us to deal with various compositions of set-valued and single-valued mappings and to establish maximal monotonicity of structurally composed operators under the validity of the corresponding qualification conditions. We refer the reader to both volumes of \cite{M1} for numerous applications of the coderivative calculus to different issues of variational analysis, optimization, and control while not related to monotonicity.

Note also that a similar full calculus is available in \cite{M1} for the {\em normal} limiting coderivative, which is defined by scheme (\ref{2.5}) with replacing the strong convergence $y\to w$ therein by the weak convergence in $X$. However, the corresponding positive-semidefiniteness condition in terms of the normal coderivative is only sufficient (together with the imposed hypomonotonicity) for the maximal monotonicity of $T$ outside finite-dimensional spaces. The proof of the necessity part given in Theorem~\ref{Coro1} does not hold true for the normal coderivative, since we cannot pass to the limit in the inequality $\la z_k,w_k\ra\ge 0$ when both sequences $\{z_k\}$ and $\{w_k\}$ converge only weakly in $X$ as $k\to\infty$.}
\end{Remark}

The following one-dimensional example shows that the hypomonotonicity conditions in (ii) in Theorem~\ref{thm1}, Theorem~\ref{thm2}, and Theorem~\ref{Coro1} are essential for the obtained coderivative characterizations of maximal monotonicity.

\begin{Example}[\bf hypomonotonicity is essential]\label{exhypo}{\rm Given $\kk\ge 0$, define the set-valued mapping $T\colon\R\tto\R$ with full domain given by:
$$
T(x):=\kk x+[0,1]\;\mbox{ for all }\;x\in\R.
$$
It is easy to calculate directly by the definitions (or using elementary calculus) that
$$
D^*T(u,v)(w)=\Hat D^*T(u,v)(w)=\begin{cases}\{0\}&\mbox{if }\;w=0,\;v-\kk u\in(0,1),\\
\{\kk w\}&\mbox{if }\;w\ge 0,\;v-ku=0,\\
\{\kk w\}&\mbox{if }\;w\le 0,\;v-ku=1,\\
\emp&\mbox{otherwise}.
\end{cases}
$$
Thus both coderivative conditions (\ref{Cod}) and (\ref{Lim-Cod}) are satisfied. However, $T$ is not monotone. The reason is that this mapping is not semilocally hypomonotone.}
\end{Example}

As consequences of the obtained results, we derive in the next corollary the corresponding regular and limiting coderivative characterizations of strong maximum monotonicity for set-valued mappings in Hilbert spaces. Recall that $T:X\tto X$ is (globally) {\em strongly maximal monotone} on $X$ with modulus $\kk>0$ if it is maximal monotone and the shifted mapping $T-\kk I$ is monotone on $X$. It follows from the classical Minty theorem that $T$ is strongly maximal monotone with modulus $\kk$ if and only if  $T-\kk I$ is maximal monotone.

\begin{Corollary}[\bf coderivative characterizations of strong maximal monotonicity]\label{Coro2} Let $T:X\tto X$ be a set-valued mapping with closed graph. The following are equivalent:

{\bf (i)} $T$ is strongly maximal monotone on $X$ with modulus $\kk>0$.

{\bf (ii)} $T$ is  hypomonotone on $X$ and for any $(u,v)\in\gph T$ we have
\begin{eqnarray*}
\la z,w\ra\ge\kk\|w\|^2\;\mbox{ whenever }\;z\in\Hat D^*T(u,v)(w),\;w\in X.
\end{eqnarray*}

{\bf (iii)} $T$ is hypomonotone on $X$ and for any $(u,v)\in\gph T$ we have
\begin{eqnarray*}
\la z,w\ra\ge \kk\|w\|^2\;\mbox{ whenever }\;z\in D_M^*T(u,v)(w),\;w\in X.
\end{eqnarray*}
If in addition the operator domain $\dom T$ is convex, then the $($global$)$ hypomonotonicity in assertions {\bf(ii)} and {\bf(iii)} can be equivalently replaced by the semilocal one.
\end{Corollary}
\noindent{\bf Proof.} Define $S:=T-\kk I$ and get from the corresponding coderivative sum rules in \cite[Theorem~1.62]{M1} the equalities
\[
\hat D^*T(u,v)(w)=\Hat D^*S(u,v-\kk u)(w)+\kk w\;\mbox{ and }\;D_M^*T(u,v)(w)=D_M^*S(u,v-\kk u)(w)+\kk w
\]
for all $(u,v)\in\gph T$ and $w\in X$. Thus the validity of (ii) (resp.\ (iii)) for $T$ is equivalent to the fulfillment of all the conditions in Theorem~\ref{thm1}(ii) (resp.\ in Theorem~\ref{Coro1}(ii)) for the operator $S$; it is obvious for hypomonotonicity. Applying now Theorem~\ref{thm1} and Theorem~\ref{Coro1}, respectively, we get that either assertion (ii) or (iii) of this corollary is equivalent to the maximal monotonicity of $S$. Since the latter is equivalent to the strong maximal monotonicity of $T$ with modulus $\kk$, we complete the proof of the corollary.\endproof

Note that a certain localized regular coderivative characterization of the {\em local} strong maximal monotonicity for set-valued mappings with respect to {\em product} neighborhoods (see the discussion right before Definition~\ref{LH}) has been recently obtained in \cite[Theorem~3.4]{MN} while being independent of the global characterizations in Corollary~\ref{Coro2}.

\section{Second-Order Characterizations of Convexity}

In this section we apply the obtained coderivative characterizations of maximal monotonicity and second-order subdifferential constructions to characterize {\em convexity} and {\em strong convexity} for the remarkable class of lower-${\cal C}^2$ functions.

Recall \cite[Definition~10.29]{rw} that a function $f\colon\R^n\rightarrow\R$ is {\em lower-${\cal C}^k$} with $k\in\IN\cup\{\infty\}$ if for each $\ox\in\R^n$ there is a neighborhood $V$ of $\bar x$ on which $\ph$ admits the representation
$$
f(x)=\max\limits_{t\in T}f_t(x),\quad x\in V,
$$
where the functions $f_t$ are of class ${\cal C}^k$ on $V$, where $T$ is compact, and where $f_t(x)$ and all their partial derivatives in $x$ through order $k$ depend continuously on $(t,x)\in T\times V$. This class of functions  introduced by Rockafellar \cite{r82} is among the favorable classes of functions in variational analysis and optimization. Many nice properties and equivalent descriptions of such functions can be found, e.g.,  in \cite{adt,r82,rw}.  As shown in \cite[Corollary~10.34]{rw}, for each $k>2$ the class of lower-${\cal C}^k$ functions coincides with the class of lower-${\cal C}^2$ functions. However, the latter is a proper subclass of lower-${\cal C}^1$ functions. In fact, between the class of lower-${\cal C}^1$ functions and the class of lower-${\cal C}^2$ functions there are the classes of lower-${\cal C}^{1,\alpha}$ functions, $0<\alpha\le 1$, which have been studied recently by Daniilidis and Malick \cite{AM25}.

The next theorem is the main result of this section, which provides complete second-order characterizations of convexity for lower-${\cal C}^2$ functions in finite dimensions.

\begin{Theorem}[\bf second-order subdifferential characterizations of convexity for the class of lower-${\cal C}^2$ functions]\label{a-thm1}  Let
$f\colon\R^n\to\R$ be a lower-${\cal C}^2$ function. Then the following assertions are equivalent:

$({\bf i})$ $f$ is convex on $\R^n$.

$({\bf ii})$ For each $(u,v)\in\gph\partial f$ we have the condition
\begin{eqnarray}\label{Lim-Conv}
\la z,w\ra\ge 0\;\mbox{ whenever }\;z\in\partial^2f(u,v)(w),\;w\in\R^n.
\end{eqnarray}

$({\bf iii})$ For each $(u,v)\in\gph\partial f$ we have the condition
\begin{eqnarray}\label{Lim-Conv1}
\la z,w\ra\ge 0\;\mbox{ whenever }\;z\in\breve\partial^2f(u,v)(w),\;w\in\R^n.
\end{eqnarray}
\end{Theorem}
\noindent{\bf Proof.} To verify implication (i)$\Longrightarrow$(ii), observe that the subdifferential operator $\partial\ph$ is maximal monotone for any convex l.s.c.\ function by the classical result of convex analysis. Hence condition (\ref{Lim-Conv}) follows from implication (i)$\Longrightarrow$(ii) in Theorem~\ref{Coro1} by construction (\ref{2.7}) of the second-order subdifferential in finite dimensions. The next implication (ii)$\Longrightarrow$(iii) of the theorem is obvious due to the inclusion
$$
\breve\partial^2f(u,v)(w)\subset\partial^2f(u,v)(w)\;\mbox{ for every }\;(u,v)\in\gph\partial f,\;w\in\R^n.
$$

To prove finally implication (iii)$\Longrightarrow$(i), suppose that $f$ is lower-${\cal C}^2$ and that condition (\ref{Lim-Conv1}) holds. It follows from \cite[Example~12.28]{rw} that the subdifferential mapping $T(x):=\partial\ph(x)$ for the lower-${\cal C}^2$ function $\ph$ is semilocally hypomonotone with $\dom T=\R^n$. Then (c) amounts to saying that for any $(u,v)\in\gph T$ we have
$$
\la z,w\ra\ge 0\;\mbox{ for all }\;z\in D^*T(u,v)(w),\;w\in\R^n.
$$
Since the set $\partial\ph$ is closed in this setting, we deduce from Theorem~\ref{thm2} that $T$ is monotone, and thus $\ph$ is convex by the result of \cite{cjt}; see also \cite[Theorem~3.56]{M1}.\endproof

\begin{Remark} [\bf discussion on second-order subdifferential characterizations of convexity]\label{other-conv} {\rm The following comments are in order:

{\bf (i)} Consider the pointwise maximum of finitely many ${\cal C}^2$ functions
\begin{eqnarray}\label{fin-max}
f(x):=\max\big\{f_1(x),\ldots,f_m(x)\big\},\quad x\in\R,
\end{eqnarray}
which surely belongs to the class of lower-${\cal C}^2$ functions. Based in the recent precise calculation \cite{eh} of the second-order subdifferential $\partial^2\ph$ of (\ref{fin-max}) via the generated functions $\ph_i$ and the appropriate index subsets of $\{1,\ldots,m\}$, we can express the second-order characterization (\ref{Lim-Conv}) entirely in terms of the initial data of (\ref{fin-max}).

{\bf (ii)} Second-order subdifferential characterizations of convexity have been recently obtained in \cite{CCYY,CH,CT,CY11,CYY14} for functions that may not be necessarily lower-${\cal C}^2$, and thus they are generally independent of Theorem~\ref{a-thm1}. Note also that it is not clear by now whether the characterizations of Theorem~\ref{a-thm1} hold if $f$ is merely lower-${\cal C}^1$ or lower-${\cal C}^{1,\alpha}$ with $0<\alpha\le 1$.}
\end{Remark}

Finally in this section, we present the second-order subdifferential characterizations of strong convexity for lower-${\cal C}^2$ functions, which can be treated as consequences of Theorem~\ref{a-thm1}. Recall that $f$ is {\em strongly convex} on $\R^n$ with modulus $\kappa>0$ if
\begin{eqnarray*}
f\big(t\lm x+(1-\lm)y\big)\le t f(x)+(1-\lm)f(y)-\frac{\kappa}{2}\lm(1-\lm)\|x-y\|^2
\end{eqnarray*}
for all $x,y\in\R$ and $\lm\in(0,1)$. It is well known that the strong convexity of $f$ is equivalent to the convexity of the shifted function \begin{eqnarray}\label{shift}
g(x):=f(x)-\frac{\kappa}{2}\|x\|^2,\quad x\in\R^n.
\end{eqnarray}

\begin{Corollary} [\bf second-order subdifferential characterizations of strong convexity for lower-${\cal C}^2$ functions]\label{strong-conv}
Let $f\colon\R^n\to\R$  be lower-${\cal C}^2$. The following are equivalent.

{\bf (i)} $f$ is strongly convex on $\R^n$ with modulus $\kappa>0$.

{\bf (ii)} We have the second-order subdifferential condition
\begin{eqnarray*}
\la z,w\ra\ge\kappa\|w\|^2\;\mbox{ for all }\;z\in\partial^2f(u,v)(w),\;(u,v)\in\gph\partial f,\;w\in\R^n.
\end{eqnarray*}

{\bf (iii)} We have the modified second-order subdifferential condition
\begin{eqnarray*}
\la z,w\ra\ge\kappa\|w\|^2\;\mbox{ for all }\;z\in\breve\partial^2f(u,v)(w),\;(u,v)\in\gph\partial f,\;w\in\R^n.
\end{eqnarray*}
\end{Corollary}
\noindent{\bf Proof.} It can be derived from Theorem~\ref{a-thm1} by applying it to the shifted function (\ref{shift}) and taking into account the obvious subdifferential relationship
$$
\partial g(x)=\partial f(x)-\kappa x,\quad x\in\R^n.
$$
On the other hand, we can justify the results by applying the characterizations of strong maximal monotonicity from Corollary~\ref{Coro2} similarly to the proof of Theorem~\ref{a-thm1}.\endproof

\section{Concluding Remarks}

The main results of this paper provide coderivative characterizations of maximal monotonicity for set-valued mappings under one of the following  assumptions:  {\bf (a)} the mapping is globally hypomonotone without imposing the convexity of its domain, and {\bf(b)} the mapping is semilocally  hypomonotone with convex domain. Example~\ref{exhypo} shows that removing hypomonotonicity may destroy these characterizations. However, it is proved in \cite{CT} that for single-valued mappings hypomonotonicity can be replaced by continuity. Thus the {\em first} open question is to clarify
what is common for both hypomonotonicity of set-valued mappings and continuity of single-valued ones. We intend to develop general results in this direction, which unify these two requirements.

The {\em second} area of the promising future research is to find an umbrella, which covers all the second-order subdifferential characterizations of convexity for functions developed in the previous investigations \cite{CCYY,CH,CT,CY11,CYY14} and those obtained in Theorem~\ref{a-thm1} for lower-${\cal C}^2$ functions. So far we can deduce \cite[Theorem~4.1]{CCYY} from Theorem~\ref{a-thm1} while the other major results of the aforementioned papers seem to be independent. It is also desired to obtain second-order subdifferential characterizations of convexity for the important classes of lower-${\cal C}^1$ and lower-${\cal C}^{1,\alpha}$ $(0<\alpha\le 1)$ functions.

The {\em third} and probably most important direction of the future research is employing the limiting coderivative calculus to derive from the pointwise coderivative characterization (\ref{Lim-Cod}) verifiable conditions for {\em preserving maximal monotonicity} under various combinations (including sums, compositions, etc.) of set-valued and single-valued maximal monotone operators. The classical result in this vein is Rockafellar's theorem \cite{r70} about maximal monotonicity of sums of maximal monotone operators under certain interiority or local boundedness assumptions. It seems that the coderivative characterizations of maximal monotonicity obtained in our paper open a new gate (in Hilbert spaces so far) to proceed in this direction via verifying qualification conditions that ensure the validity of the corresponding coderivative calculus rules from \cite{M1}.

\end{document}